\documentclass{gtart}

\def\ifplaintex{\expandafter\ifx\csname documentclass\endcsname\relax}

\def\gtp{{\mathsurround=0pt\it $\cal G\mskip-2mu$eometry \&\ 
$\cal T\!\!$opology $\cal P\!$ublications}}  

\def\recd{{\small Received:\qua\receiveddate\ifx\reviseddate\relax
\else\qquad Revised:\qua\reviseddate\fi\par}} 

\def\lognumber#1{\def\thelognumber{#1}}
\def\volumenumber#1{\def\thevolumenumber{#1}}
\def\volumeyear#1{\def\thevolumeyear{#1}}
\def\papernumber#1{\def\thepapernumber{#1}}
\def\pagenumbers#1#2{\def\startpage{#1}\def\finishpage{#2}}
\def\published#1{\def\publishdate{#1}}

\def\received#1{\def\receiveddate{#1}}

\def\accepted#1{\def\accepteddate{#1}}
\def\asciititle#1{\def\theasciititle{#1}}
\def\covertitle#1{\def\thecovertitle{#1}}

\long\def\asciiabstract#1{\long\def\theasciiabstract{#1}}

\let\\\par\let\thelognumber\relax\let\thevolumenumber\relax
\let\thepapernumber\relax\let\thevolumeyear\relax\let\startpage\relax
\let\finishpage\relax\let\publishdate\relax\let\receiveddate\relax
\let\reviseddate\relax\let\accepteddate\relax\let\theasciititle\relax
\let\thecovertitle\relax\let\theasciiauthors\relax
\let\theasciiabstract\relax

\let\theasciiemail\relax

\ifplaintex
\font\logobig=cmssbx10 scaled 3836
\font\logomed=cmssbx10 scaled 2557
\else
\font\logobig=cmssbx10 scaled 4200
\font\logomed=cmssbx10 scaled 2800
\fi

\long\def\makeagttitle{   
\count0=\startpage
\agt\hfill      
\hbox to 45truept{\vbox to 0pt{\vglue -13truept{\logomed A\kern -.37em{\logobig 
T}\kern -.38em G}\vss}\hss}
\break
{\small Volume \thevolumenumber\ (\thevolumeyear)
\startpage--\finishpage\nl
Published: \publishdate}

\vglue .25truein

{\parskip=0pt\leftskip 0pt plus
1fil\def\\{\par\smallskip}{\Large\bf\thetitle}\par\medskip} \vglue
0.05truein

{\parskip=0pt\leftskip 0pt plus 1fil\def\\{\par}{\sc\theauthors}
\par\medskip}%
 
\vglue 0.03truein 

{\small\leftskip 25truept\rightskip 25truept{\bf Abstract}\stdspace\theabstract

{\bf AMS Classification}\stdspace\theprimaryclass
\ifx\thesecondaryclass\relax\else; \thesecondaryclass\fi\par
{\bf Keywords}\stdspace \thekeywords\par}\vglue 7truept

}   

\ifplaintex
\font\phead=cmsl9 scaled 950
\font\pnum=cmbx10 scaled 913
\font\pfoot=cmsl9 scaled 950
\headline{\vbox to 0pt{\vskip -4.5mm\line{\small\phead\ifnum
\count0=\startpage ISSN 1472-2739 (on-line) 1472-2747 (printed)
\hfill {\pnum\folio}\else\ifodd\count0\def\\{ }%
\ifx\theshorttitle\relax\thetitle\else\theshorttitle\fi\hfill{\pnum\folio}
\else\def\\{ and }{\pnum\folio}\hfill\ifx\theshortauthors\relax\theauthors
\else\theshortauthors\fi\fi\fi}\vss}}
\footline{\vbox to 0pt{\vglue 0mm\line{\small\pfoot\ifnum\count0=\startpage
\copyright\ \gtp\hfill\else
\agt, Volume \thevolumenumber\ (\thevolumeyear)\hfill\fi}\vss}}
\else
\font=cmsl9 scaled 1050
\font\lnum=cmbx10 
\font=cmsl9 scaled 1050
\makeatletter
\def\@oddhead{{\small\lhead\ifnum\count0=\startpage ISSN 1472-2739 
(on-line) 1472-2747 (printed)\hfill {\lnum\number\count0}\else\ifodd\count0
\def\\{ }\ifx\theshorttitle\relax \thetitle \else\theshorttitle\fi\hfill
{\lnum\number\count0}\else\def\\{ and }{\lnum\number\count0}
\hfill\ifx\theshortauthors\relax 
\theauthors\else\theshortauthors\fi\fi\fi}}\def\@evenhead{\@oddhead}
\def\@oddfoot{\small\lfoot\ifnum\count0=\startpage\copyright\ \gtp\hfill\else
\agt, Volume \thevolumenumber\ (\thevolumeyear)\hfill\fi}
\def\@evenfoot{\@oddfoot}
\makeatother
\fi
\let\maketitlepage\makeagttitle
\let\makeshorttitle\maketitlepage
\let\maketitle\maketitlepage

\newwrite\gtoutfile
\long\gdef\makeheadfile{  
{\def\\{, }\def\s{ }
\immediate\openout\gtoutfile head.xxx
\immediate\write\gtoutfile{To: math@arxiv.org}
\immediate\write\gtoutfile{Subject: put OR rep NNNNN:ppppp}
\immediate\write\gtoutfile{--text follows this line--}
\immediate\write\gtoutfile{Proxy-for: \ifx\theasciiauthors\relax
\theauthors\else\theasciiauthors\fi\s<\ifx\theasciiemail\relax\theemail\else\theasciiemail\fi>}
\immediate\write\gtoutfile{\noexpand\\}
\immediate\write\gtoutfile{Authors: \ifx\theasciiauthors\relax
\theauthors\else\theasciiauthors\fi}
{\def\\{ }\immediate\write\gtoutfile{Title: \ifx\theasciititle\relax
\thetitle\else\theasciititle\fi}}
\immediate\write\gtoutfile{Subj-class: GT or SG, GR etc}
\immediate\write\gtoutfile{MSC-class: \theprimaryclass\ifx\thesecondaryclass\relax\else, \thesecondaryclass\fi}
\immediate\write\gtoutfile{Journal-ref: Algebr. Geom. Topol. \thevolumenumber\s
(\thevolumeyear) \startpage-\finishpage}
\immediate\write\gtoutfile{Comments: Published by Algebraic and
Geometric Topology at}
\immediate\write\gtoutfile{\s\s\s  http://www.maths.warwick.ac.uk/agt/AGTVol\thevolumenumber/agt-\thevolumenumber-\thepapernumber.abs.html}
\immediate\write\gtoutfile{\noexpand\\}
\immediate\write\gtoutfile{}
\ifx\theasciiabstract\relax
\immediate\write\gtoutfile{\theabstract}\else
\immediate\write\gtoutfile{\theasciiabstract}\fi
\immediate\write\gtoutfile{}
\immediate\write\gtoutfile{\noexpand\\}
\immediate\write\gtoutfile{}
\immediate\closeout\gtoutfile}}  

\def\maketitlepage{\makeagttitle\makeheadfile}
\let\makeshorttitle\maketitlepage
\let\maketitle\maketitlepage

\def\ifplaintex{\expandafter\ifx\csname documentclass\endcsname\relax}

\def\gtp{{\mathsurround=0pt\it $\cal G\mskip-2mu$eometry \&\ 
$\cal T\!\!$opology $\cal P\!$ublications}}  

\def\recd{{\small Received:\qua\receiveddate\ifx\reviseddate\relax
\else\qquad Revised:\qua\reviseddate\fi\par}} 

\def\lognumber#1{\def\thelognumber{#1}}
\def\volumenumber#1{\def\thevolumenumber{#1}}
\def\volumeyear#1{\def\thevolumeyear{#1}}
\def\papernumber#1{\def\thepapernumber{#1}}
\def\pagenumbers#1#2{\def\startpage{#1}\def\finishpage{#2}}
\def\published#1{\def\publishdate{#1}}

\def\received#1{\def\receiveddate{#1}}

\def\accepted#1{\def\accepteddate{#1}}
\def\asciititle#1{\def\theasciititle{#1}}
\def\covertitle#1{\def\thecovertitle{#1}}

\long\def\asciiabstract#1{\long\def\theasciiabstract{#1}}

\let\\\par\let\thelognumber\relax\let\thevolumenumber\relax
\let\thepapernumber\relax\let\thevolumeyear\relax\let\startpage\relax
\let\finishpage\relax\let\publishdate\relax\let\receiveddate\relax
\let\reviseddate\relax\let\accepteddate\relax\let\theasciititle\relax
\let\thecovertitle\relax\let\theasciiauthors\relax
\let\theasciiabstract\relax

\let\theasciiemail\relax

\ifplaintex
\font\logobig=cmssbx10 scaled 3836
\font\logomed=cmssbx10 scaled 2557
\else
\font\logobig=cmssbx10 scaled 4200
\font\logomed=cmssbx10 scaled 2800
\fi

\long\def\makeagttitle{   
\count0=\startpage
\agt\hfill      
\hbox to 45truept{\vbox to 0pt{\vglue -13truept{\logomed A\kern -.37em{\logobig 
T}\kern -.38em G}\vss}\hss}
\break
{\small Volume \thevolumenumber\ (\thevolumeyear)
\startpage--\finishpage\nl
Published: \publishdate}

\vglue .25truein

{\parskip=0pt\leftskip 0pt plus
1fil\def\\{\par\smallskip}{\Large\bf\thetitle}\par\medskip} \vglue
0.05truein

{\parskip=0pt\leftskip 0pt plus 1fil\def\\{\par}{\sc\theauthors}
\par\medskip}%
 
\vglue 0.03truein 

{\small\leftskip 25truept\rightskip 25truept{\bf Abstract}\stdspace\theabstract

{\bf AMS Classification}\stdspace\theprimaryclass
\ifx\thesecondaryclass\relax\else; \thesecondaryclass\fi\par
{\bf Keywords}\stdspace \thekeywords\par}\vglue 7truept

}   

\ifplaintex
\font\phead=cmsl9 scaled 950
\font\pnum=cmbx10 scaled 913
\font\pfoot=cmsl9 scaled 950
\headline{\vbox to 0pt{\vskip -4.5mm\line{\small\phead\ifnum
\count0=\startpage ISSN 1472-2739 (on-line) 1472-2747 (printed)
\hfill {\pnum\folio}\else\ifodd\count0\def\\{ }%
\ifx\theshorttitle\relax\thetitle\else\theshorttitle\fi\hfill{\pnum\folio}
\else\def\\{ and }{\pnum\folio}\hfill\ifx\theshortauthors\relax\theauthors
\else\theshortauthors\fi\fi\fi}\vss}}
\footline{\vbox to 0pt{\vglue 0mm\line{\small\pfoot\ifnum\count0=\startpage
\copyright\ \gtp\hfill\else
\agt, Volume \thevolumenumber\ (\thevolumeyear)\hfill\fi}\vss}}
\else
\font=cmsl9 scaled 1050
\font\lnum=cmbx10 
\font=cmsl9 scaled 1050
\makeatletter
\def\@oddhead{{\small\lhead\ifnum\count0=\startpage ISSN 1472-2739 
(on-line) 1472-2747 (printed)\hfill {\lnum\number\count0}\else\ifodd\count0
\def\\{ }\ifx\theshorttitle\relax \thetitle \else\theshorttitle\fi\hfill
{\lnum\number\count0}\else\def\\{ and }{\lnum\number\count0}
\hfill\ifx\theshortauthors\relax 
\theauthors\else\theshortauthors\fi\fi\fi}}\def\@evenhead{\@oddhead}
\def\@oddfoot{\small\lfoot\ifnum\count0=\startpage\copyright\ \gtp\hfill\else
\agt, Volume \thevolumenumber\ (\thevolumeyear)\hfill\fi}
\def\@evenfoot{\@oddfoot}
\makeatother
\fi
\let\maketitlepage\makeagttitle
\let\makeshorttitle\maketitlepage
\let\maketitle\maketitlepage

\newwrite\gtoutfile
\long\gdef\makeheadfile{  
{\def\\{, }\def\s{ }
\immediate\openout\gtoutfile head.xxx
\immediate\write\gtoutfile{To: math@arxiv.org}
\immediate\write\gtoutfile{Subject: put OR rep NNNNN:ppppp}
\immediate\write\gtoutfile{--text follows this line--}
\immediate\write\gtoutfile{Proxy-for: \ifx\theasciiauthors\relax
\theauthors\else\theasciiauthors\fi\s<\ifx\theasciiemail\relax\theemail\else\theasciiemail\fi>}
\immediate\write\gtoutfile{\noexpand\\}
\immediate\write\gtoutfile{Authors: \ifx\theasciiauthors\relax
\theauthors\else\theasciiauthors\fi}
{\def\\{ }\immediate\write\gtoutfile{Title: \ifx\theasciititle\relax
\thetitle\else\theasciititle\fi}}
\immediate\write\gtoutfile{Subj-class: GT or SG, GR etc}
\immediate\write\gtoutfile{MSC-class: \theprimaryclass\ifx\thesecondaryclass\relax\else, \thesecondaryclass\fi}
\immediate\write\gtoutfile{Journal-ref: Algebr. Geom. Topol. \thevolumenumber\s
(\thevolumeyear) \startpage-\finishpage}
\immediate\write\gtoutfile{Comments: Published by Algebraic and
Geometric Topology at}
\immediate\write\gtoutfile{\s\s\s  http://www.maths.warwick.ac.uk/agt/AGTVol\thevolumenumber/agt-\thevolumenumber-\thepapernumber.abs.html}
\immediate\write\gtoutfile{\noexpand\\}
\immediate\write\gtoutfile{}
\ifx\theasciiabstract\relax
\immediate\write\gtoutfile{\theabstract}\else
\immediate\write\gtoutfile{\theasciiabstract}\fi
\immediate\write\gtoutfile{}
\immediate\write\gtoutfile{\noexpand\\}
\immediate\write\gtoutfile{}
\immediate\closeout\gtoutfile}}  

\def\maketitlepage{\makeagttitle\makeheadfile}
\let\makeshorttitle\maketitlepage
\let\maketitle\maketitlepage

\def\ifplaintex{\expandafter\ifx\csname documentclass\endcsname\relax}

\def\gtp{{\mathsurround=0pt\it $\cal G\mskip-2mu$eometry \&\ 
$\cal T\!\!$opology $\cal P\!$ublications}}  

\def\recd{{\small Received:\qua\receiveddate\ifx\reviseddate\relax
\else\qquad Revised:\qua\reviseddate\fi\par}} 

\def\lognumber#1{\def\thelognumber{#1}}
\def\volumenumber#1{\def\thevolumenumber{#1}}
\def\volumeyear#1{\def\thevolumeyear{#1}}
\def\papernumber#1{\def\thepapernumber{#1}}
\def\pagenumbers#1#2{\def\startpage{#1}\def\finishpage{#2}}
\def\published#1{\def\publishdate{#1}}

\def\received#1{\def\receiveddate{#1}}

\def\accepted#1{\def\accepteddate{#1}}
\def\asciititle#1{\def\theasciititle{#1}}
\def\covertitle#1{\def\thecovertitle{#1}}

\long\def\asciiabstract#1{\long\def\theasciiabstract{#1}}

\let\\\par\let\thelognumber\relax\let\thevolumenumber\relax
\let\thepapernumber\relax\let\thevolumeyear\relax\let\startpage\relax
\let\finishpage\relax\let\publishdate\relax\let\receiveddate\relax
\let\reviseddate\relax\let\accepteddate\relax\let\theasciititle\relax
\let\thecovertitle\relax\let\theasciiauthors\relax
\let\theasciiabstract\relax

\let\theasciiemail\relax

\ifplaintex
\font\logobig=cmssbx10 scaled 3836
\font\logomed=cmssbx10 scaled 2557
\else
\font\logobig=cmssbx10 scaled 4200
\font\logomed=cmssbx10 scaled 2800
\fi

\long\def\makeagttitle{   
\count0=\startpage
\agt\hfill      
\hbox to 45truept{\vbox to 0pt{\vglue -13truept{\logomed A\kern -.37em{\logobig 
T}\kern -.38em G}\vss}\hss}
\break
{\small Volume \thevolumenumber\ (\thevolumeyear)
\startpage--\finishpage\nl
Published: \publishdate}

\vglue .25truein

{\parskip=0pt\leftskip 0pt plus
1fil\def\\{\par\smallskip}{\Large\bf\thetitle}\par\medskip} \vglue
0.05truein

{\parskip=0pt\leftskip 0pt plus 1fil\def\\{\par}{\sc\theauthors}
\par\medskip}%
 
\vglue 0.03truein 

{\small\leftskip 25truept\rightskip 25truept{\bf Abstract}\stdspace\theabstract

{\bf AMS Classification}\stdspace\theprimaryclass
\ifx\thesecondaryclass\relax\else; \thesecondaryclass\fi\par
{\bf Keywords}\stdspace \thekeywords\par}\vglue 7truept

}   

\ifplaintex
\font\phead=cmsl9 scaled 950
\font\pnum=cmbx10 scaled 913
\font\pfoot=cmsl9 scaled 950
\headline{\vbox to 0pt{\vskip -4.5mm\line{\small\phead\ifnum
\count0=\startpage ISSN 1472-2739 (on-line) 1472-2747 (printed)
\hfill {\pnum\folio}\else\ifodd\count0\def\\{ }%
\ifx\theshorttitle\relax\thetitle\else\theshorttitle\fi\hfill{\pnum\folio}
\else\def\\{ and }{\pnum\folio}\hfill\ifx\theshortauthors\relax\theauthors
\else\theshortauthors\fi\fi\fi}\vss}}
\footline{\vbox to 0pt{\vglue 0mm\line{\small\pfoot\ifnum\count0=\startpage
\copyright\ \gtp\hfill\else
\agt, Volume \thevolumenumber\ (\thevolumeyear)\hfill\fi}\vss}}
\else
\font=cmsl9 scaled 1050
\font\lnum=cmbx10 
\font=cmsl9 scaled 1050
\makeatletter
\def\@oddhead{{\small\lhead\ifnum\count0=\startpage ISSN 1472-2739 
(on-line) 1472-2747 (printed)\hfill {\lnum\number\count0}\else\ifodd\count0
\def\\{ }\ifx\theshorttitle\relax \thetitle \else\theshorttitle\fi\hfill
{\lnum\number\count0}\else\def\\{ and }{\lnum\number\count0}
\hfill\ifx\theshortauthors\relax 
\theauthors\else\theshortauthors\fi\fi\fi}}\def\@evenhead{\@oddhead}
\def\@oddfoot{\small\lfoot\ifnum\count0=\startpage\copyright\ \gtp\hfill\else
\agt, Volume \thevolumenumber\ (\thevolumeyear)\hfill\fi}
\def\@evenfoot{\@oddfoot}
\makeatother
\fi
\let\maketitlepage\makeagttitle
\let\makeshorttitle\maketitlepage
\let\maketitle\maketitlepage

\newwrite\gtoutfile
\long\gdef\makeheadfile{  
{\def\\{, }\def\s{ }
\immediate\openout\gtoutfile head.xxx
\immediate\write\gtoutfile{To: math@arxiv.org}
\immediate\write\gtoutfile{Subject: put OR rep NNNNN:ppppp}
\immediate\write\gtoutfile{--text follows this line--}
\immediate\write\gtoutfile{Proxy-for: \ifx\theasciiauthors\relax
\theauthors\else\theasciiauthors\fi\s<\ifx\theasciiemail\relax\theemail\else\theasciiemail\fi>}
\immediate\write\gtoutfile{\noexpand\\}
\immediate\write\gtoutfile{Authors: \ifx\theasciiauthors\relax
\theauthors\else\theasciiauthors\fi}
{\def\\{ }\immediate\write\gtoutfile{Title: \ifx\theasciititle\relax
\thetitle\else\theasciititle\fi}}
\immediate\write\gtoutfile{Subj-class: GT or SG, GR etc}
\immediate\write\gtoutfile{MSC-class: \theprimaryclass\ifx\thesecondaryclass\relax\else, \thesecondaryclass\fi}
\immediate\write\gtoutfile{Journal-ref: Algebr. Geom. Topol. \thevolumenumber\s
(\thevolumeyear) \startpage-\finishpage}
\immediate\write\gtoutfile{Comments: Published by Algebraic and
Geometric Topology at}
\immediate\write\gtoutfile{\s\s\s  http://www.maths.warwick.ac.uk/agt/AGTVol\thevolumenumber/agt-\thevolumenumber-\thepapernumber.abs.html}
\immediate\write\gtoutfile{\noexpand\\}
\immediate\write\gtoutfile{}
\ifx\theasciiabstract\relax
\immediate\write\gtoutfile{\theabstract}\else
\immediate\write\gtoutfile{\theasciiabstract}\fi
\immediate\write\gtoutfile{}
\immediate\write\gtoutfile{\noexpand\\}
\immediate\write\gtoutfile{}
\immediate\closeout\gtoutfile}}  

\def\maketitlepage{\makeagttitle\makeheadfile}
\let\makeshorttitle\maketitlepage
\let\maketitle\maketitlepage

\def\ifplaintex{\expandafter\ifx\csname documentclass\endcsname\relax}

\def\gtp{{\mathsurround=0pt\it $\cal G\mskip-2mu$eometry \&\ 
$\cal T\!\!$opology $\cal P\!$ublications}}  

\def\recd{{\small Received:\qua\receiveddate\ifx\reviseddate\relax
\else\qquad Revised:\qua\reviseddate\fi\par}} 

\def\lognumber#1{\def\thelognumber{#1}}
\def\volumenumber#1{\def\thevolumenumber{#1}}
\def\volumeyear#1{\def\thevolumeyear{#1}}
\def\papernumber#1{\def\thepapernumber{#1}}
\def\pagenumbers#1#2{\def\startpage{#1}\def\finishpage{#2}}
\def\published#1{\def\publishdate{#1}}

\def\received#1{\def\receiveddate{#1}}

\def\accepted#1{\def\accepteddate{#1}}
\def\asciititle#1{\def\theasciititle{#1}}
\def\covertitle#1{\def\thecovertitle{#1}}

\long\def\asciiabstract#1{\long\def\theasciiabstract{#1}}

\let\\\par\let\thelognumber\relax\let\thevolumenumber\relax
\let\thepapernumber\relax\let\thevolumeyear\relax\let\startpage\relax
\let\finishpage\relax\let\publishdate\relax\let\receiveddate\relax
\let\reviseddate\relax\let\accepteddate\relax\let\theasciititle\relax
\let\thecovertitle\relax\let\theasciiauthors\relax
\let\theasciiabstract\relax

\let\theasciiemail\relax

\ifplaintex
\font\logobig=cmssbx10 scaled 3836
\font\logomed=cmssbx10 scaled 2557
\else
\font\logobig=cmssbx10 scaled 4200
\font\logomed=cmssbx10 scaled 2800
\fi

\long\def\makeagttitle{   
\count0=\startpage
\agt\hfill      
\hbox to 45truept{\vbox to 0pt{\vglue -13truept{\logomed A\kern -.37em{\logobig 
T}\kern -.38em G}\vss}\hss}
\break
{\small Volume \thevolumenumber\ (\thevolumeyear)
\startpage--\finishpage\nl
Published: \publishdate}

\vglue .25truein

{\parskip=0pt\leftskip 0pt plus
1fil\def\\{\par\smallskip}{\Large\bf\thetitle}\par\medskip} \vglue
0.05truein

{\parskip=0pt\leftskip 0pt plus 1fil\def\\{\par}{\sc\theauthors}
\par\medskip}%
 
\vglue 0.03truein 

{\small\leftskip 25truept\rightskip 25truept{\bf Abstract}\stdspace\theabstract

{\bf AMS Classification}\stdspace\theprimaryclass
\ifx\thesecondaryclass\relax\else; \thesecondaryclass\fi\par
{\bf Keywords}\stdspace \thekeywords\par}\vglue 7truept

}   

\ifplaintex
\font\phead=cmsl9 scaled 950
\font\pnum=cmbx10 scaled 913
\font\pfoot=cmsl9 scaled 950
\headline{\vbox to 0pt{\vskip -4.5mm\line{\small\phead\ifnum
\count0=\startpage ISSN 1472-2739 (on-line) 1472-2747 (printed)
\hfill {\pnum\folio}\else\ifodd\count0\def\\{ }%
\ifx\theshorttitle\relax\thetitle\else\theshorttitle\fi\hfill{\pnum\folio}
\else\def\\{ and }{\pnum\folio}\hfill\ifx\theshortauthors\relax\theauthors
\else\theshortauthors\fi\fi\fi}\vss}}
\footline{\vbox to 0pt{\vglue 0mm\line{\small\pfoot\ifnum\count0=\startpage
\copyright\ \gtp\hfill\else
\agt, Volume \thevolumenumber\ (\thevolumeyear)\hfill\fi}\vss}}
\else
\font=cmsl9 scaled 1050
\font\lnum=cmbx10 
\font=cmsl9 scaled 1050
\makeatletter
\def\@oddhead{{\small\lhead\ifnum\count0=\startpage ISSN 1472-2739 
(on-line) 1472-2747 (printed)\hfill {\lnum\number\count0}\else\ifodd\count0
\def\\{ }\ifx\theshorttitle\relax \thetitle \else\theshorttitle\fi\hfill
{\lnum\number\count0}\else\def\\{ and }{\lnum\number\count0}
\hfill\ifx\theshortauthors\relax 
\theauthors\else\theshortauthors\fi\fi\fi}}\def\@evenhead{\@oddhead}
\def\@oddfoot{\small\lfoot\ifnum\count0=\startpage\copyright\ \gtp\hfill\else
\agt, Volume \thevolumenumber\ (\thevolumeyear)\hfill\fi}
\def\@evenfoot{\@oddfoot}
\makeatother
\fi
\let\maketitlepage\makeagttitle
\let\makeshorttitle\maketitlepage
\let\maketitle\maketitlepage

\newwrite\gtoutfile
\long\gdef\makeheadfile{  
{\def\\{, }\def\s{ }
\immediate\openout\gtoutfile head.xxx
\immediate\write\gtoutfile{To: math@arxiv.org}
\immediate\write\gtoutfile{Subject: put OR rep NNNNN:ppppp}
\immediate\write\gtoutfile{--text follows this line--}
\immediate\write\gtoutfile{Proxy-for: \ifx\theasciiauthors\relax
\theauthors\else\theasciiauthors\fi\s<\ifx\theasciiemail\relax\theemail\else\theasciiemail\fi>}
\immediate\write\gtoutfile{\noexpand\\}
\immediate\write\gtoutfile{Authors: \ifx\theasciiauthors\relax
\theauthors\else\theasciiauthors\fi}
{\def\\{ }\immediate\write\gtoutfile{Title: \ifx\theasciititle\relax
\thetitle\else\theasciititle\fi}}
\immediate\write\gtoutfile{Subj-class: GT or SG, GR etc}
\immediate\write\gtoutfile{MSC-class: \theprimaryclass\ifx\thesecondaryclass\relax\else, \thesecondaryclass\fi}
\immediate\write\gtoutfile{Journal-ref: Algebr. Geom. Topol. \thevolumenumber\s
(\thevolumeyear) \startpage-\finishpage}
\immediate\write\gtoutfile{Comments: Published by Algebraic and
Geometric Topology at}
\immediate\write\gtoutfile{\s\s\s  http://www.maths.warwick.ac.uk/agt/AGTVol\thevolumenumber/agt-\thevolumenumber-\thepapernumber.abs.html}
\immediate\write\gtoutfile{\noexpand\\}
\immediate\write\gtoutfile{}
\ifx\theasciiabstract\relax
\immediate\write\gtoutfile{\theabstract}\else
\immediate\write\gtoutfile{\theasciiabstract}\fi
\immediate\write\gtoutfile{}
\immediate\write\gtoutfile{\noexpand\\}
\immediate\write\gtoutfile{}
\immediate\closeout\gtoutfile}}  

\def\maketitlepage{\makeagttitle\makeheadfile}
\let\makeshorttitle\maketitlepage
\let\maketitle\maketitlepage

\lognumber{36}
\volumenumber{2}
\volumeyear{2002}
\papernumber{36}
\pagenumbers{897}{919}
\received{4 February 2002}
\accepted{21 August 2002}
\published{20 October 2002}

\usepackage{amssymb,amsmath}
\input xy
\xyoption{all}

\newtheorem{Thm}{Theorem}[section]
\newtheorem{Cor}[Thm]{Corollary}
\newtheorem{Prop}[Thm]{Proposition}
\newtheorem{Lem}[Thm]{Lemma}
\theoremstyle{remark}
\newtheorem{Dfn}[Thm]{Definition}
\newtheorem{Rmk}[Thm]{Remark}

\newcommand{\fix}{\mbox{\rm Fix } \phi}
\newcommand{\fp}[1]{\mbox{\rm Fix } #1}

\begin{document}

\title[Maximal index automorphisms]{Maximal index automorphisms of
free groups with\\\vglue-7pt\\no attracting fixed points on the
boundary\\are Dehn twists}
\covertitle{Maximal index automorphisms of
free groups with\\no attracting fixed points on the boundary\\are Dehn
twists}
\asciititle{Maximal index automorphisms of
free groups with no attracting fixed points on the boundary are Dehn
twists}

\author{Armando Martino}

\address{Department of Mathematics, University College Cork\\Cork, Ireland}

\email{A.Martino@ucc.ie}

\begin{abstract}
In this paper we define a quantity called the {\em rank} of an
outer automorphism of a free group which is the same as the index
introduced in \cite{GJLL} without the count of fixed points on the
boundary. We proceed to analyze outer automorphisms of maximal
rank and obtain results analogous to those in \cite{CT1}. We also
deduce that all such outer automorphisms can be represented by
Dehn twists, thus proving the converse to a result in \cite{CL},
and indicate a solution to the conjugacy problem when such
automorphisms are given in terms of images of a basis, thus
providing a moderate extension to the main theorem of \cite{CL} by
somewhat different methods.
\end{abstract}

\asciiabstract{ In this paper we define a quantity called the rank of
an outer automorphism of a free group which is the same as the index
introduced in [D. Gaboriau, A. Jaeger, G. Levitt and M. Lustig, `An
index for counting fixed points for automorphisms of free groups',
Duke Math. J. 93 (1998) 425-452] without the count of fixed points on
the boundary. We proceed to analyze outer automorphisms of maximal
rank and obtain results analogous to those in [D.J. Collins and
E. Turner, `An automorphism of a free group of finite rank with
maximal rank fixed point subgroup fixes a primitive element', J. Pure
and Applied Algebra 88 (1993) 43-49]. We also deduce that all such
outer automorphisms can be represented by Dehn twists, thus proving
the converse to a result in [M.M. Cohen and M. Lustig, `The conjugacy
problem for Dehn twist automorphisms of free groups', Comment
Math. Helv.  74 (1999) 179-200], and indicate a solution to the
conjugacy problem when such automorphisms are given in terms of images
of a basis, thus providing a moderate extension to the main theorem of
Cohen and Lustig by somewhat different methods.}

\primaryclass{20E05, 20E36}
\keywords{Free group, automorphism}

\maketitle

\section{Introduction}

The celebrated result of \cite{BH} showed that for any
automorphism of a finitely generated free group the rank of its
fixed subgroup is at most that of the rank of the ambient free
group. In \cite{CT2} a detailed analysis and description was
obtained for those automorphisms whose fixed subgroup has the
largest possible rank -- maximal rank automorphisms.

The paper of \cite{CL} introduced a class of automorphisms called
Dehn Twists (defined below) and showed that these have maximal
index with no attracting fixed points on the boundary. In that
work, the conjugacy problem for Dehn Twists is also solved and it
is shown, by using the results of \cite{CT2}, that a maximal rank
automorphism can be represented by a Dehn Twist.

In this paper we define a notion of rank for outer automorphisms
which generalises the notion of rank of the fixed subgroup. (In
fact this notion of rank is implicit in \cite{BH}.) Alternatively,
this rank can be thought of as the index of an outer automorphism,
as described in \cite{GJLL}, but with the change that the fixed
points on the boundary are not counted.

We then proceed to generalise the results of \cite{CT2} to the
class of maximal rank outer automorphisms and obtain a normal form
similar to the one obtained there. Moreover, we show that any such
outer automorphism can be realised as a Dehn Twist. Thus the class
of Dehn Twists and that of maximal rank outer automorphisms
coincide.

In \cite{AM} the normal form of \cite{CT2} is used to provide a
solution to the conjugacy problem for maximal rank automorphisms.
We also observe that since the normal form we obtain is so similar
to that in \cite{CT2}, it is possible to use the same proof
 to provide a solution to the conjugacy problem for Dehn
Twists by entirely different means to those of \cite{CL}.
Moreover, this solution would take as input data a Dehn Twist
described purely in terms of its action on a basis rather than by
graph of groups data as required in \cite{CL} hence giving an
extension to that result.

\section{Preliminaries}

\subsection{Outer automorphisms and index}

Throughout $F_n$ shall denote the free group of rank $n$. Here the
rank is the minimal number of generators and is the same as the
number of free generators. The rank of a subgroup, $H$, of $F_n$
is the least cardinality of the generating sets of the subgroup
and is denoted $r(H)$.

$Aut(F_n)$ is the group of automorphisms of $F_n$. $Inn(F_n)$ will
be the subgroup of inner automorphisms and $Out(F_n) =
Aut(F_n)/Inn(F_n)$ the group of outer automorphisms of $F_n$. We
use the notation $\gamma_g$ to denote conjugation by $g$. Thus $w
\gamma_g = g^{-1} w g$ for all $w \in F_n$. (We will write
automorphisms on the right, although the topological
representatives below will be written on the left).

We shall think of an outer automorphism $\Phi $ of $F_n$ as a
coset, and as such it will be a set of automorphisms any two of
which differ by conjugation by some element.

A similarity class in $\Phi$ will be an equivalence class under
the equivalence relation, $\phi \sim \psi $ if and only if $\phi =
\gamma_g \psi \gamma_{g^{-1}}$ for some $g \in F_n$. Note that we
think of this as an equivalence relation on $\Phi$ where two
equivalent automorphisms are called similar.

It is important to note that, unlike the situation with matrices,
two automorphisms are {\em not} called similar if they are
conjugate. They are only similar if and only if they are conjugate
by an {\em inner} automorphism. This follows the terminology of
\cite{GJLL}. In \cite{CL}, the same concept is denoted by the
phrase {\em conjugate up to inner automorphisms} and some authors
use instead the term {\em isogredience}.

The fixed subgroup of an automorphism $\phi \in Aut(F_n)$ is the
subgroup $\fix = \{ w \in F_n : w  \phi =w \}$ and by \cite{BH}
has rank at most $n$. Given an outer automorphism $\Phi $, one can
find (infinitely many) representatives $\phi_i$ of the distinct
similarity classes in $\Phi $.
 It is clear that similar automorphisms have fixed subgroups
 of the same rank thus
the following (possibly infinite) quantity is well defined.

\begin{Dfn}
The rank of an outer automorphism $\Phi$ of $F_n$ is the sum
$$
r(\Phi):=1+\sum max(0,r(\fp{\phi_i})-1),$$  where the sum is taken
over a set $\{\phi_i\}$ of representatives, one for each
similarity class of $\Phi$.
\end{Dfn}

The following Theorem is proved in \cite{GLL} and is in fact
equivalent to the main Theorem of \cite{BH}.
\begin{Thm}\label{max}{\rm\cite{GLL}}\qua
For every $\Phi \in Out(F_n)$, $r(\Phi) \leq n$.
\end{Thm}
 An immediate observation is that only finitely many of the
$\phi_i$ have fixed subgroup of rank greater than one.

This observation leads us to the following definition.
\begin{Dfn}
\upshape For any $\Phi \in Out(F_n)$, let $s(\Phi)$ denote the
number of distinct similarity classes in $\Phi $ which have fixed
subgroup of rank at least $2$. By Theorem \ref{max}, this quantity
is always finite.
\end{Dfn}

In \cite{GJLL} the {\em index} of an outer automorphism $i(\Phi)$
is defined. This has the same definition as the rank of $\Phi$
defined above with the addition of a term which counts attracting
fixed points on the boundary of $F_n$ . Thus it is clear that
$r(\Phi) \leq i(\Phi)$. In \cite{GJLL} it is shown that $i(\Phi)
\leq n$ for every $\Phi \in Out(F_n)$. \footnote{In fact the index
described in \cite{GJLL} is one less than the one we refer to. The
change is merely to emphasise the parallels between Theorems in
$Out(F_n)$ and $Aut(F_n)$.}
 Also in \cite{CL} it is
shown that if $\Phi$ is represented by a Dehn twist automorphism
(definitions below), then $i(\Phi )=n$ and $\Phi $ has no
attracting points on the boundary. This is the same as saying that
$r(\Phi)=n$. In this paper we prove the converse of this result.
Namely that if $r(\Phi)=n$ then $\Phi$ can be represented by a
Dehn twist automorphism.

In fact the main Theorem of \cite{BH}, proved a conjecture of
Scott's who formulated it after proving the following:
\begin{Thm}\label{scott}{\rm\cite{DS}}\qua
If an automorphism $\phi$ of $Aut(F_n)$ has finite
order, then $\fix $ is a free factor of $F_n$.
\end{Thm}

\begin{Cor}
\label{frf} For any $\phi \in Aut(F_n)$ and any integer $m \geq
1$, $\fix $ is a free factor of $\fp{\phi^m}$.
\end{Cor}
\proof By Theorem \ref{max}, $\fp{\phi^m}$ is of finite rank and as
$\phi$ and $\phi^m$ commute, $\phi$ leaves $\fp{\phi^m}$
invariant. Since it acts as a finite order automorphism, the
result follows from Theorem \ref{scott}.\endproof

This will have important consequences for us. The construction
used in the proof of the following Proposition is due to G.
Levitt.
\begin{Lem}
\label{finite} Consider an outer automorphism $\Phi \in Out(F_n)$
of finite order. If $\Phi$ is not the identity then $r(\Phi)<n$.
\end{Lem}
\proof We may find a maximal set of automorphisms $\phi_1,\ldots
,\phi_k \in \Phi$ in distinct similarity classes all of whose
fixed subgroups have rank at least $2$. (Note that $k=s(\Phi)$.)
Thus $r(\Phi)=1+ \sum_{j=1}^k (r(\fp{\phi_j})-1)$.

Also, we know that $\Phi $ has order $m$ in $Out(F_n)$ for some $m
\geq 2$. Hence every automorphism in $\Phi^m$ is inner and by
Theorem \ref{scott} this implies that ${\phi_j}^m=1 \in Aut(F_n)$
for all $1 \leq j \leq k$.

Pick a basis $x_1,\ldots ,x_n$ for $F_n$ and find $g_2,\ldots,g_k$
so that $\phi_j \gamma_{g_j} =\phi_1$. (By definition the $\phi_j$
differ by inner automorphisms.)

Consider a free group of rank $n+k-1$, $F$, with basis $x_1,\ldots
,x_n,\ldots ,x_{n+k-1}$ where we identify $F_n$ with $\langle
x_1,\ldots ,x_n \rangle$. Define an automorphism $\phi$ of $F$ by
setting $\phi|_{F_n}=\phi_1$ and $(x_{n+j-1}) \phi =x_{n+j-1}
g_j$, $2 \leq j \leq k$.

First note that $\phi $ cannot fix certain words. $\phi $ cannot
fix any word of the form $x_j w {x_{j'}}^{-1}$ for $j \not=j' \geq
n+1$ and $w \in F_n$, for if it did then this would imply that
$$
\begin{array}{rcl}
g_j ( w )\phi_1 {g_{j'}}^{-1} & = & w \ \Rightarrow \\
w^{-1}  (w) \phi_j  & = & g_{j'} {g_j}^{-1} \ ,\mbox{\rm since }
 \phi_j \gamma_{g_j}=\phi_1. \\
\end{array}
$$
This last equality is not possible since it would mean that
$$
\begin{array}{rcl}
{\gamma_w}^{-1} \phi_j \gamma_w & = & \phi_j \gamma_{
(w^{-1})\phi_j w} \\
& = & \phi_j \gamma_{g_j {g_{j'}}^{-1}} \\
& = & \phi_j \gamma_{g_j} \gamma_{{g_{j'}}^{-1}} \\
& = & \phi_1 \gamma_{{g_{j'}}^{-1}} \\
& = & \phi_{j'},
\end{array}
$$
and by construction these automorphisms are not similar.

Also, $\phi $ cannot fix any word of the form $x_j w$, for $j \geq
n+1$ and $w \in F_n$ for then we get $g_j = w (w^{-1})\phi$. This
would imply the similarity of $\phi_1$ and $\phi_j$ as $\gamma_w
\phi_1 {\gamma_w}^{-1} = \phi_j$ again reaching a contradiction.

We are now in a position to determine $\fix$.

\medskip

\noindent {\bf Claim}
$$\fix= \fp{\phi_1} \ast_{j=2}^k x_{n+j-1} (\fp{\phi_j}) {x_{n+j-1}}^{-1}.$$ It is
clear that the term on the right hand side is a subgroup of
$\fix$.
 Consider a word $w
\not\in F_n$, of shortest length fixed by $\phi$ and not of the
form given above. We can write such a word as,
$$
w_0 {x_{j_1}}^{\epsilon_1} w_1 {x_{j_2}}^{\epsilon_2} w_2 \ldots
{x_{j_p}}^{\epsilon_p} w_p
$$
where each $j_i  \geq n+1$, $\epsilon_i = \pm 1$ and $w_i \in F_n$
and $p \geq 1$. We proceed by a simple cancellation argument.

If $\epsilon_1=-1$ it is easy to see that $x_{j_1} {w_0}^{-1}$
must be fixed, giving a contradiction as above. Hence
$\epsilon_1=1$ and since this implies that $w_0$ is fixed we may
assume that $w_0=1$ by minimality of the length of $w$.

 Now, if $\epsilon=1$ and either $p=1$ or
$\epsilon_2=1$, then $x_{j_1} w_1$ must be fixed, again a
contradiction. Thus $p$ must be at least $2$, $\epsilon_1=1$ and
$\epsilon_2=-1$ leading us to the conclusion that $x_{j_1} w_1
{x_{j_2}}^{-1}$  is fixed. The only way that $\phi $ can fix a
word of this type is if $j_1=j_2$ and $w_1 \in \fp{\phi_{j_1}}$.
This contradicts the minimality of $w$ and proves the claim.

Hence,
$$ r(\fix) = \sum_{j=1}^k r(\fp{\phi_j}) = r(\Phi) +k-1.$$
However, each $\phi_j$ has order $m$ and since $\phi_j= \phi_1
\gamma_{{g_j}^{-1}}$, this implies that $(g_j^{-1})\phi_1^{m-1}
(g_j^{-1}) \phi_1^{m-2} \ldots g_j^{-1}=1$.

Taking inverses we get that $g_j (g_j) \phi_1 (g_j) \phi_1^2\ldots
(g_j)\phi_1^{m-1}=1$ and hence that $x_j \phi^m = x_j$. Since
$\phi^m|_{F_n}=\phi_1^m=1$ we know that $\phi^m=1$. Clearly, $\phi
\not= 1$ and so by Theorem \ref{scott},  $r(\fix) < r(F) = n+k-1$.
Thus $r(\Phi) = r(\fix) -k+1 < n $, completing the proof of the
Lemma. \endproof

Consider an outer automorphism $\Phi$. Clearly if $\phi, \psi \in
\Phi $ are similar, then $\phi^m, \psi^m \in \Phi^m$ will also be
similar. The converse on the other hand need not be true. However,
if we concentrate on those similarity classes with non-trivial
fixed subgroup, the next proposition tells us that the only way
these similarity classes can get `collapsed' in a power is if the
rank of the outer automorphism increases.
\begin{Prop}
\label{count} Consider $\Phi \in Out(F_n)$ and let $\phi,\psi \in
\Phi$ be non-similar automorphisms each with non-trivial fixed
subgroup. If for some integer $m$, $\phi^m$ and $\psi^m$ are
similar, then $r(\Phi) < r(\Phi^m)$.
\end{Prop}
\proof We first choose a collection of automorphisms $\phi_1,\ldots
,\phi_k \in \Phi$ in distinct similarity classes each with
non-trivial subgroup so that each of $\phi $ and $\psi$ is similar
to some automorphism on the list. Additionally, we enlarge the
list so that any automorphism in $\Phi$ which has fixed subgroup
of rank at least $2$ is similar to one on the list.

 After a rearrangement we may find integers $k_1,\ldots ,k_s$ so
that $\phi^m_i$ is similar to $\phi_j^m$ if and only if $ k_p \leq
i,j < k_{p+1}$ for some $1 \leq k_p < k_s$. In other words, we
list representatives of similarity classes for $\Phi$ so that only
consecutive similarity classes get collapsed in $\Phi^m$. As a
consequence, the automorphisms $\phi_{k_1}^m,\phi_{k_2}^m,\ldots
,\phi_{k_s}^m$ form a set of representatives of distinct
similarity classes of $\Phi^m$ with non-trivial fixed subgroup.
Thus
$$
r(\Phi^m) \geq 1+\sum_{p=1}^s (r(\fp{\phi_{k_p}^m})-1).$$
By changing representatives for similarity classes in $\Phi$ we
may in fact assume that $\phi_i^m=\phi_j^m$ whenever $k_p \leq i,j
< k_{p+1}$. Thus if $k_p \leq j < k_{p+1}$, then $\fp{\phi_j}$ is
a subgroup of $\fp{\phi_j^m}=\fp{\phi_{k_p}^m}:=H_p$. In fact, as
$\phi_j$ and $\phi_j^m$ commute, $H_p$ is invariant under $\phi_j$
which restricts to a finite order automorphism on this subgroup.

Also, if we write $\phi_j= \phi_{k_p} \gamma_g $ then,
$$
\begin{array}{rcl}
 \phi_j^m & = & \phi_{k_p}^m \gamma_{(g)\phi_{k_p}^{m-1}
 (g)\phi_{k_p}^{m-2}\ldots (g)\phi_{k_p} g}  \\
 & = & \phi_j^m
\gamma_{(g)\phi_{k_p}^{m-1}
 (g)\phi_{k_p}^{m-2}\ldots(g)\phi_{k_p} g} \\
\end{array}
$$
and hence $(g)\phi_{k_p}^{m-1}
 (g)\phi_{k_p}^{m-2}\ldots (g)\phi_{k_p} g=1$. Looking at the image of
 this element under $\phi_{k_p}$ we note that $
 (g)\phi_{k_p}^m g^{-1}=1$. In other words, $g \in \fp{\phi_{k_p}^m}=H_p$ and
 the two automorphisms in question induce the same outer
 automorphism when restricted to $H_p$.

Note that if $g$ were to be fixed by $\phi_{k_p}$ then this would
imply that that $g=1$ and so that $\phi_j=\phi_{k_p}$. Hence if
$k_{p+1}-k_p>1$ then $H_p=\fp{\phi_{k_p}^m}$ is strictly larger
than $\fp{\phi_{k_p}}$ and in particular, the outer automorphism
induced by $\phi_{k_p}$ on $H_p$ cannot be the identity.

Thus let $\Psi_p$ be the outer automorphism induced by the
restriction of $\phi_{k_p}$ to $H_p$. By the comments above, we
know that $\Psi_p$ is a finite order outer automorphism and that,
$$
r(\Psi_p) \geq 1 + \sum_{j=k_p}^{k_{p+1} -1} (r(\fp{\phi_j})-1).$$
Note that $r(\Psi_p)$ is always bounded above by $r(H_p)$ by
Theorem \ref{max}, however if the number of terms in the sum on
the right hand side is greater than one we know that $\Psi_p \not=
1$ and hence we may apply Lemma \ref{finite} to deduce that
$r(\Psi_p) < r(H_p)$. In fact, our hypothesis guarantees that for
some $p$ this will be the case, and so
$$
1 + \sum_{p=1}^s ( r(\Psi_p) - 1) < 1 + \sum_{i=1}^s( r(H_p)-1).
$$
As the left hand side is bounded below by $r(\Phi)$ and the right
hand side is bounded above by $r(\Phi^m)$ this concludes the
proof. \endproof

Recall that $w \in F_n$ is $\phi $ periodic if $(w)\phi^m=w$ for
some $m \not=0$.
\begin{Cor}
\label{noper} Let $\Phi \in Out(F_n)$ have maximal rank and
suppose that $w \in F_n$ is $\phi $ periodic for some $\phi \in
\Phi$ with non-trivial fixed subgroup. Then $w \in \fix$.
\end{Cor}
\proof Choose an integer $m$ such that $w \phi^m=w$. Since
$r(\Phi)=n$ we deduce, by Proposition \ref{count},  that if $\phi,
\psi \in \Phi$ have non-trivial fixed subgroup, then they are
similar if and only if $\phi^m,\psi^m \in \Phi^m$ are similar.
 In particular this
means that if $r(\fix)<r(\fp{\phi^m})$ then $r(\Phi)<r(\Phi^m)$, a
contradiction. But by Corollary \ref{frf}, if
$r(\fix)=r(\fp{\phi^m})$ then in fact $\fix = \fp{\phi^m}$. This
completes the proof. \endproof



\section{Relative train track maps}

We use the relative train track maps of Bestvina and Handel to
analyze an outer automorphism. We recap some of the properties of
relative train track maps.

A relative train track map is a self homotopy equivalence, $f$, of
a graph $G$ which maps vertices to vertices and edges to paths.
Such an $f$ is called a topological representative.

Note that if the image $f(e)$ of an edge $e$ is {\em not}
homotopic to a trivial path relative endpoints then we may replace
$f$ with a homotopically equivalent map that is locally injective
on the interior of $e$. The process of doing this for each edge is
called {\em tightening} $f$. Relative train track maps are always
tightened.

The graph $G$ has no valence one vertices and is maximally
filtered in the sense that it has subgraphs $$ \emptyset =G_0
\subseteq G_1 \subseteq \cdots \subseteq G_m=G, $$ where each
$G_i$ is an $f$-invariant subgraph and if $f(G_r) \not\subseteq
G_{r-1}$ then there is no $f$-invariant subgraph strictly between
$G_{r-1}$ and $G_r$. The closure of $G_r \backslash G_{r-1}$ is
denoted by $H_r$ and is called the $r^{th}$ stratum.

On labelling the edges of the $r^{th}$ stratum, $e_1,\ldots ,e_k$,
one can form the $r^{th}$ transition matrix, $M_r$, whose $(i,j)$
entry is the number of times that that $f(e_i)$ crosses $e_j$ (in
either direction).

If $f(G_r) \subseteq G_{r-1}$ then $M_r$ is a zero matrix and
$H_r$ is called a zero stratum. If $M_r$ is a permutation matrix,
then $H_r$ is called a level stratum. Otherwise $H_r$ is called an
exponential stratum.

\begin{Rmk}
\upshape In order for $f,G$ to be a relative train track map
further conditions need to be imposed on the exponential strata.
However, we shall only need to consider, by Proposition
\ref{noexp}, those maps with no exponential strata, in which case
relative train track maps are precisely those topological
representatives which are tight and maximally filtered.
\end{Rmk}

A Nielsen path (NP) is a path in $G$ which is fixed by $f$ up to
homotopy relative endpoints. An indivisible Nielsen path (INP) is
an NP which cannot be written (non-trivially) as the concatenation
of NP's. In \cite{BH} it is shown that every NP can be written
uniquely as a product of INP's.

A path in $G$ is said to have height $r$ if it is contained in
$G_r$ but not $G_{r-1}$. In \cite{BH} it is shown that there is at
most one INP of height $r$ for each $r$. Note that this uses the
property known as {\em stability} in \cite{BH}. We shall also
assume throughout that isolated fixed points of $f$ are in fact
vertices of $G$. The following remarks are part of the analysis of
\cite{BH} and in particular, the proof of Proposition~6.3.

\begin{Rmk}(\cite{BH}, pp48-49)\qua
\label{level} \upshape If there is an INP, $\rho $, of height $r$,
then $H_r$ cannot be a zero stratum. Furthermore, if $H_r$ is a
level stratum then it must consist of a single edge $E$, with
$f(E)=E u $ for some path $u$ in $G_{r-1}$. In that case, $\rho $
must be of the form $E \beta$ or $E \beta \bar{E}$ for some path
$\beta $ in $G_{r-1}$.
\end{Rmk}

A graph $\Sigma$ along with a map $p\co  \Sigma \to G$ may then be
constructed such that,

\begin{enumerate}
\item $p$ maps vertices of $\Sigma$ to vertices of $G$ fixed by
$f$,
\item $p$ maps edges of $\Sigma $ to INPs, and
\item every NP in $G$ is the image (under $p$) of a path in
$\Sigma$.
\end{enumerate}

In fact, $\Sigma$ is constructed so that its vertices can be
regarded as being precisely those vertices of $G$ which are fixed
by $f$. Also, following \cite{BH}, one may define certain
subgraphs of $\Sigma$.
\begin{Dfn}
\label{stratsub} Let $\Sigma_r$ to be the (not necessarily
connected) maximal subgraph of $\Sigma$ which maps to $G_r$ under
$p$.
\end{Dfn}
\begin{Dfn}
\label{nielsensub} For every vertex $v$ of $G$, which is fixed by
$f$, $\Sigma^v$ is the component of $\Sigma$ containing $v$.
\end{Dfn}

\begin{Rmk}(\cite{BH}, p48)
\label{single} \upshape The graph $\Sigma_r$ differs from
$\Sigma_{r-1}$ by at most a single edge when there is an INP of
height $r$.
\end{Rmk}

For a connected graph, $G$, define the rank of $G$, $r(G)$ to be
the rank of $\pi_1(G)$. For an arbitrary graph define the {\em
reduced rank} to be $$ \widetilde{r}(G)=1+\sum max(0,r(G_k)-1)$$
where the sum ranges over the components of $G$.

 It is shown in \cite{BH}, p48,
that

(1)\qua $\widetilde{r}(\Sigma) \leq \widetilde{r}(G)=r(G)$ and, 

(2)\qua $\widetilde{r}(\Sigma_r) \leq \widetilde{r}(G_r)$.

By Remark \ref{single} above we also have that 

(3)\qua $\widetilde{r}(\Sigma_{r+1}) \leq \widetilde{r}(\Sigma_r)+1$.

\begin{Rmk}(\cite{BH}, p48)\qua
\label{inc} \upshape
 Note that if $\widetilde{r}(\Sigma_{r+1}) =
\widetilde{r}(\Sigma_r)+1$ then there is an edge in $\Sigma$
(which maps to an INP, $\rho_r$ of height $r$ in $G$) and which
has endpoints in (possibly one) non-contractible components of
$\Sigma_{r}$. Thus $\rho_r$ also has endpoints in (possibly one)
non-contractible components of $G_r$ and so
$\widetilde{r}(G_{r+1}) \geq \widetilde{r}(G_r)+1$.
\end{Rmk}

\subsection*{Representing automorphisms}

Let $f$ be a relative train track map on the graph $G$. Suppose
that $v$ is a vertex of $G$ and $\mu $ a path in $G$ from $f(v)$
to $v$. Then $\pi_1(f,\mu)$ will denote the induced isomorphism of
$\pi_1(G,v)$ that sends the closed path $\alpha$ at $v$ to
$\bar{\mu}  f(\alpha) \mu$. Here $\bar{\mu} $ denotes the inverse
path to $\mu$. Write $\pi_1(f,v)$ in the case where $v$ is fixed
by $f$ and $\mu$ is the trivial path at $v$.

Let $R_n$ denote the graph with one vertex, $*$, and $n$ edges,
called the rose and identify $F_n$ with $\pi_1(R_n,*)$. We say
that an outer automorphism $\Phi \in Out(F_n)$ is represented by
the relative train track map $f$ on $G$ if there is a homotopy
equivalence, $\tau\co R_n \to G$ such that the following diagram
commutes up to free homotopy:
$$ \xymatrix{ R_n \ar[r]^\tau \ar[d]^\Phi &  G \ar[d]^f  \\
R_n \ar[r]^\tau & G
}
$$
Note that we are identifying $\Phi$ with a self homotopy
equivalence of $R_n$.

Given a representation of $\Phi $ as above, we say that $\phi \in
\Phi$ is point represented at $v$ (by $f$, $G$, $\tau$) if $v$ is
fixed by $f$ and there is a path, $\alpha $, from $\tau(*)$ to $v$
such that the following diagram commutes,
$$
\xymatrix{ \pi_1(R_n,*) \ar[r]^{\pi_1(\tau,\alpha)} \ar[d]^\phi &
\pi_1(G,v) \ar[d]^{\pi_1(f,v)}  \\
\pi_1(R_n,*) \ar[r]^{\pi_1(\tau,\alpha)} &  \pi_1(G,v) }
$$
where $\pi_1(\tau,\alpha)$ is induced by the map which sends the
path $g \subset R_n$ to $\bar{\alpha} \tau(g) \alpha$.

 It is shown in \cite{BH} that every outer automorphism,
$\Phi$, is represented by a relative train track map, $f,G,\tau$.
Furthermore we have:
\begin{Prop}[Corollary 2.2, \cite{BH}]
\label{point} If an (ordinary) automorphism, $\phi \in \Phi$ has
fixed subgroup of rank at least $2$, then this automorphism will
be point represented by $(f,G,\tau)$. \end{Prop}

Note that if $\phi $ is point represented at $v$, then any
automorphism similar to $\phi $ will also be point represented at
$v$. Also, if there is a Nielsen path between the vertices $v$
and $v'$, then $\phi $ will be point represented at $v'$.
Conversely, suppose that $\phi$ is point represented at both $v$
and $v'$, with paths $\alpha, \alpha'$ from $\tau(*)$ to $v,v'$
respectively, then $\bar{\alpha} \alpha'$ is a Nielsen path from
$v$ to $v'$.

\begin{Rmk}
\label{vert} \upshape Hence, bringing this together, if $\phi,
\phi' \in \Phi $ are both point represented at $v$ and $v'$
respectively (by $f$,$G$,$\tau$) then they are similar if and only
if there is a Nielsen path from $v$ to $v'$. In the case where all
NP's are closed, each fixed vertex will determine a distinct
similarity class of $\Phi$.
\end{Rmk}

If $\Phi $ is represented by $f,G,\tau$ and $\phi \in \Phi$ is
point represented at $v$ then (Definition~\ref{nielsensub})
$r(\Sigma^v)=r(\fix)$ and $\widetilde{r} (\Sigma) = r(\Phi)$. In
fact the map $p\co \Sigma \to G$ induces isomorphisms from
$\pi_1(\Sigma^v)$ to $\fix$ whenever $\phi $ is point represented
at $v$.

In the case where $\Phi$ has maximal rank (and is represented by
$f,G,\tau$) then $\widetilde{r}(\Sigma)=r(G)$. By Remarks
\ref{single} and \ref{inc} we deduce that:
\begin{Lem}
\label{equals} If $\Phi $ has maximal rank then,
$$\widetilde{r}(\Sigma_k)=\widetilde{r}(G_k) \ \mbox{\rm for all }  k.$$
\end{Lem}

\section{Good Representatives}

From now on $\Phi $ will be a maximal rank outer automorphism of
$F_n$. We shall show in this section that every such outer
automorphism has a relative train track map representative with
good properties. The first step is to observe:

\begin{Prop}[Prop 4, \cite{CT2}]
\label{noexp} If $\Phi$ has maximal rank then any relative train
track map representative has no exponential strata.
\end{Prop}

In fact the Proposition in \cite{CT2} relates to automorphisms
(not outer!) of $F_n$ with fixed subgroup of rank $n$. However,
the only hypothesis used is that
$\widetilde{r}(\Sigma_k)=\widetilde{r}(G_k)$ for all $k$, and
hence the proof there applies equally in our situation.

In order to find a good relative train track map representative,
we need to perform a certain operation as follows. Suppose that
$f\co G \to G$ is a relative train track map and that $H_k=\{E\}$,
where $f(E)=E u$, and $u$ is a path in $G_{k-1}$. For any path
$\alpha$ in $G_{k-1}$ with initial vertex the same as the terminal
vertex of $E$ we can define a new graph $G'$ by replacing $E$ with
an edge $E'$ with the same initial vertex as $v$ and whose
terminal vertex is the same as that of $\alpha$. Every edge of
$G-\{E\}$ is naturally identified with that of $G'-\{E'\}$. We can
then define $f'\co G' \to G'$ so as to agree with $f$ (up to
homotopy) on $G-\{E\}$, and so that $f'(E') \simeq f(E)
\bar{\alpha} u f( \alpha)$ and $f'$ is tight.  The homotopy
equivalence $p\co G \to G'$ which is the `identity' on $G-\{E\}$
and sends $E$ to $E' \bar{\alpha}$ gives the following commuting
diagram. $$ \xymatrix{ G \ar[r]^p \ar[d]^f & G' \ar[d]^{f'} \\ G
\ar[r]^p & G'
\\ } $$
Moreover, if we set $G'_j=p(G_j)$ then $f'$ is a relative train
track map with stratum $H'_j$ of the same type (zero, level or
exponential) as $H_j$. This operation is called {\em sliding} in
\cite{BFH} and a proof of the above statements is contained in
\cite{BFH}, Lemma 5.4.1 and is a slight variation of the
construction that appears in \cite{CT2}, Proposition 3.

Our first application of sliding is in fact precisely analogous to
that in \cite{CT2}.

\begin{Prop}
\label{goodinp}Let $\Phi \in Out(F_n)$, $n \geq 2$, have maximal
rank. Then there is a relative train track map representative,
$f,G$, for $\Phi$ in which every indivisible Nielsen path,
$\rho_k$, of height $k$ is either of the form $E \beta \bar{E}$
for some path $\beta $ in $G_{k-1}$ or $\rho_k=E$ and $E$ is a
closed loop.
\end{Prop}
\proof By Proposition \ref{noexp} and Remark \ref{level}, we only
need to consider the case where a stratum $H_k$ consists of a
single edge $E$, $f(E)=Eu$, and $E \alpha$ is an INP, for some
$u$, $\alpha$ subpaths in $G_{k-1}$. (This requires subdivision at
isolated fixed points).
 Sliding $E$ along $\alpha $
we obtain a relative train track map representative $f', G'$,
where $f(E') = E'$. If we do this in all possible cases and then
collapse any fixed edges which are not loops, we end up deleting
some strata, but otherwise still with a relative train track map
representative. It is clear that for this map, every INP is of one
of the above types. \endproof

An examination of the above proof gives a way of starting from a
representative of $\Phi $ and getting another where we have better
control of INP's. We want to have an easy way of insuring this
condition. For that we need the following: 

\begin{Dfn} \label{min} \upshape Let
$f,G$ represent the maximal rank outer automorphism $\Phi $. We
say that $f,G$ has {\em minimal complexity} if $G$ has the minimal
number of vertices amongst all representatives of $\Phi $ subject
to the restriction that $f,G$ satisfies the conclusion to
Proposition \ref{goodinp} and that all isolated fixed points are
vertices.
\end{Dfn}

\begin{Lem}
\label{goodmin} Any representative $f,G$ of $\Phi $ of minimal
complexity has the minimal number of vertices amongst all
representatives of $\Phi$ with vertices at isolated fixed points.
\end{Lem}
\proof An examination of the proof to Proposition \ref{goodinp}
shows that if a representative with vertices at isolated fixed
points does not satisfy the proposition then we perform a sliding
operation followed by the collapse of an invariant forest. Since
this cannot increase the number of vertices of the underlying
graph and cannot introduce any new isolated fixed points, we are
done. \endproof

We shall henceforth assume that our maximal rank outer
automorphism is represented by a relative train track map which
satisfies the conclusions of Proposition \ref{goodinp}.

\begin{Rmk}
\label{top} \upshape
 Suppose that $G$ has exactly $r$ strata, so that $G=G_r$.
Then since $G$ has no valence one vertices, $\widetilde{r} G_r >
\widetilde{r} G_{r-1}$. Thus by Lemma \ref{equals}, there is an
INP of height $r$. Hence there is a single edge $E$ so that
$H_r=\{E\}$ and $f(E)= E u $ with $u$ a path in $G_{r-1}$,
possibly trivial. Denote the initial vertex of $E$ by $v$ and the
terminal vertex by $w$. Let $C_1$ denote the component of
$G_{r-1}$ containing $v$ and $C_2$ the component containing $w$.
Thus if $E$ is non-separating, $C_1=C_2$. Clearly, $f(C_i)
\subseteq C_i$ and in fact it must restrict to a homotopy
equivalence in each case.
\end{Rmk}

Using this notation we can show:
\begin{Prop}
\label{rank1} Let $f,G$ be a representative of $\Phi \in Out(F_n)$
($n \geq 2$) of minimal complexity. Then the following hold:

{\rm(1)}\qua If $E$ is separating in $G$ then $C_2$ has rank at least $2$.

{\rm(2)}\qua If $C_1 $ is a rank $1$ graph, then it consists of a single
closed fixed edge and a single vertex.
\end{Prop}

\proof We start with property $1$ where we need to show
 that if $E$ is separating then
$C_2$ must have rank at least $2$. If this is not the case then
$C_2$ will have rank one and since there is an INP $E \beta
\bar{E}$, $f|_{E \cup C_2}$ is homotopic to the identity map
relative to $v$. Thus there is a map $f'$ on $G$ which also
represents $\Phi$ and which is the identity on $E \cup C_2$ ($f'$
agrees with $f$ on $C_1$). It is clear that $f'$ is also a
relative train track map (with the same strata as $f$) and has
vertices at isolated fixed points. Note that $E$ is a separating
edge fixed by $f'$. By collapsing $E$ we contradict Lemma
\ref{goodmin}.

To prove $2$, note that since $ \widetilde{r}\Sigma_r >
\widetilde{r}\Sigma_{r-1}$, $\Sigma^v_{r-1}$ must have rank at
least $1$. Hence there is always a closed Nielsen path at the
vertex $v$ contained in $C_1$. So if the rank of $C_1$ is $1$ then
we can apply the argument as above to show that, without loss of
generality, $f$ restricts to the identity on $C_1$. The only way
that this does not contradict the minimal complexity hypothesis is
if $C_1$ consists of a single fixed edge. \endproof

As an immediate corollary we get:

\begin{Cor}
\label{rank2} Let $\Phi \in Out(F_2)$ have maximal rank and $f,G$
be any relative train track map representative of minimal
complexity. Then exactly one vertex $v$ and two edges, $a,b$,
where without loss of generality $f(a) \simeq a $ and $f(b) \simeq
b a^r$ for some integer $r$.
\end{Cor}

\begin{Prop}
\label{inv} Let $f,G$ be a relative train track map representing a
maximal rank outer automorphism. Suppose that $C$ is a component
of some $G_k$ with $r(C) \geq 1$. Then $f(C) \subseteq C$ and $f$
induces a maximal rank outer automorphism on $C$.
\end{Prop}
\proof The proof is by induction on $r(G)$. If $r(G)=1$ then $C$
can only be equal to $G$ and we are done.

Consider the edge $E$ as in Remark \ref{top}. If $E$ is a
separating edge, then since it is the content of the highest
stratum, we can write $\Sigma_{r-1}$ as a disjoint union of
graphs, $\Sigma^1$ and $\Sigma^2$, where $p(\Sigma^i) \subseteq
C_i$. In other words, $\Sigma^i$ contains all the edges of
$\Sigma_{r-1}$ that map to INP's of $C_i$. It is clear by the
properties of $\Sigma $ that the rank of the outer automorphism
induced by $f|_{C_i}$ is exactly $\widetilde{r} \Sigma^i$ and
hence by Theorem \ref{max}, $\widetilde{r} \Sigma^i \leq
\widetilde{r} C_i$. However, $$
\begin{array}{rcl}
\widetilde{r}G & = & \widetilde{r} C_1 + \widetilde{r}C_2 +1 \\ &
\geq  & \widetilde{r} \Sigma^1 + \widetilde{r} \Sigma^2 +1 \\ & =
& \widetilde{r} \Sigma_{r-1} + 1 \\ & = & \widetilde{r} \Sigma \\&
= & \widetilde{r} G.\\
\end{array}
$$ The upshot of this is that each $f|_{C_i}$ induces a maximal
rank outer automorphism on $\pi_1(C_i)$. A similar argument
applies when $E$ is non-separating, where
$G-\{E\}=G_{r-1}=C_1=C_2$ to get that $f|_{C_1}$ induces a maximal
rank outer automorphism.

Now each $C_i$ inherits a filtration from $G$, namely, $C_i \cup
G_m$ is an invariant subgraph of $C_i$. Thus as any component $C$
of some $G_k$ is actually a component of $C_i \cap G_k$ (except
$G$ itself) we may use our inductive hypothesis to finish the
proof. \endproof

\begin{Thm}
\label{conj} Let $f,G$ be a relative train track map of minimal
complexity, with $r(G) \geq 2$, representing a maximal rank outer
automorphism. Suppose that for some closed path $\alpha$,
$f(\alpha) \simeq \bar{\mu} \alpha \mu$. Then there is a path
$\eta$ in $G$ such that

{\rm(1)}\qua $[\bar{\eta} \alpha \eta] \in
\pi_1(G,v)$. 

{\rm(2)}\qua $\pi_1(f,v)$ has fixed subgroup of rank at
least $2$.
\end{Thm}
\proof Let $G=G_r$ and use the notation of Remark \ref{top}.
 Note that whether or not $E$ is
separating, by Proposition \ref{goodinp},  $\Sigma$ differs from
$\Sigma_{r-1}$ by a single closed loop at the vertex $v$. Since
$\widetilde{r}\Sigma_r > \widetilde{r}\Sigma_{r-1}$, we must have
that $r(\Sigma^v_{r-1}) \geq 1$ and $r(\Sigma^v) \geq 2$.

We will first prove the Theorem in the case where $\alpha $ is not
freely homotopic to a path in $G_{r-1}$. By possibly replacing
$\alpha $ with its inverse, we may choose an $\eta$ so that
$\bar{\eta} \alpha \eta $ is a path that starts with $E$ and does
not end with $\bar{E}$. Since this is a closed path at $v$ and
$r(\Sigma^v) \geq 2$, it will suffice to show that this path is
fixed up to homotopy.

Notice that in fact every positive $f$ iterate of this path also
begins with $E$ and does not end with $\bar{E}$. Our key
observation here is that, with respect to some basis, this path
and its iterates are cyclically reduced words. If $E$ is
non-separating then choose any maximal tree that does not include
$E$. The induced basis certainly ensures that each
$[f^k(\bar{\eta} \alpha \eta)]$ are cyclically reduced.

If, on the other hand, $E$ is separating then note that $f$
induces automorphisms of $H=\pi_1(C_1,v)$ and on $K=\pi_1(C_2 \cup
\{E\},v)$. If we choose a basis for $\pi_1(G,v)$ that extends
bases for $H$ and $K$ it is then easy to see that $[\bar{\eta}
\alpha \eta] $ starts with a non-trivial word from $K$ and ends
with a non-trivial word from $H$ and the same will be true of all
its iterates. Hence in either case $[\bar{\eta} \alpha \eta] $ and
all its iterates are cyclically reduced elements of $\pi_1(G,v)$.
But since any element of a free group has only finitely cyclically
reduced conjugates, this means that $[\bar{\eta} \alpha \eta ]$ is
$\pi_1(f,v)$ periodic and hence by Corollary \ref{noper}, fixed.

This leaves us with the case where $\alpha $ is freely homotopic
to a path in $G_{r-1}=C_1 \cup C_2$.
 Consider first the case where $\alpha $ is freely homotopic to a
 path in $C_i$ and
$r(C_i) =1$. By Proposition \ref{rank1}, $C_i=C_1$ which consists
of a fixed edge loop. Thus we may choose $\eta$ so that $
\bar{\eta} \alpha \eta $ is fact a power of the fixed edge loop
and we are done in this case.

If on the other hand $\alpha \subseteq C_2$ then $E$ must be a
separating edge and the INP of height $r$ is of the form $E \beta
\bar{E}$. We can then choose $\eta $ so that $\bar{\eta} \alpha
\eta$ is a loop at $v$ homotopic to a power of that INP. (Recall
we are assuming that $r(C_2)=1$.) Again we would be done.

We finish the argument by induction on $r(G)$. If $r(G)=2$ then
either $\alpha $ is not freely homotopic to a path in $G_{r-1}$ or
$\alpha \subseteq C_i$ where $r(C_i)=1$. The arguments above deal
with each situation.

So suppose that the proposition is true for all rank less than
$r(G)$ and at least $2$. Again, if $\alpha $ is not freely
homotopic to a path in $G_{r-1}$ we are done. Hence, without loss
of generality $\alpha \subseteq C_i$ and we can assume that
$r(C_i) \geq 2$. By Proposition \ref{inv} we may apply our
induction hypothesis to complete the proof. \endproof

As an immediate consequence of the above we get:
\begin{Cor}
\label{conjalg} Let $\Phi \in Out(F_n)$, $n \geq 2$, be a maximal
rank outer automorphism fixing a conjugacy class. Then there is a
$\phi \in \Phi$ with fixed subgroup of rank at least $2$ fixing an
element of that conjugacy class.
\end{Cor}

We are now ready to show that a maximal rank outer automorphism of
$F_n$ has a representative with very good properties, analagous to
those in \cite{CT2}.

\begin{Thm}
\label{goodrep} Let $\Phi $ be a maximal rank outer automorphism
of $F_n$. Then there is a relative train track map representative
of
minimal complexity for $\Phi $, $(f,G)$, such that, 

\begin{enumerate}
\item every vertex of $G$ is fixed,
\item for every vertex, $v$
of $G$, $\pi_1(f,v)$ has fixed subgroup of rank at least $2$, and
\item (up to orientation) every edge $E$ of $G$ satisfies, $f(E)=E
\beta$, where $\beta $ is a closed Nielsen path contained in
strata lower than $E$.
\end{enumerate}
\end{Thm}

\proof
The case $n=2$ follows from Corollary \ref{rank2}, so we proceed
by induction. Our hypothesis will actually be that given any
representative of $\Phi $ of minimal complexity, a sequence of
sliding operations will transform the representative into one
which satisfies the conclusion of the Theorem.

 Start with a relative train track map
representative $f,G$ of minimal complexity and top edge $E$ as in
Remark \ref{top} so that the initial vertex of $E$ is $v$ and the
terminal vertex is $w$. We already know that the fixed subgroup of
$\pi_1(f,v)$ has rank at least $2$ and that $v$ is a fixed vertex.
In fact, we make make the following claim:

\medskip

{\bf Claim}\qua{\sl After a sliding operation we obtain  a relative train
track map representative of minimal complexity with the following
properties (notation from Remark \ref{top}):

{\rm(1)}\qua Either (i) $E$
is a fixed edge loop, or (ii) the INP of height $r$ is $E \beta
\bar{E}$, where $\beta $ is a closed Nielsen path at $w$. 

{\rm(2)}\qua In case (ii), $f(E)=E \beta^k$ for some $k$.

{\rm(3)}\qua $w$ is a fixed vertex and the fixed subgroup of
$\pi_1(f,w)$ has rank at least $2$.}

\proof[Proof of claim] The claim is immediate if $E$ is a
fixed edge loop. It also follows immediately from Proposition
\ref{rank1} if $E$ is non-separating and $r(C_1)=1$. Therefore (by
another application of \ref{rank1}) we may assume that the INP of
height $r$ is $E \beta \bar{E}$. So $\beta$ is a loop at $w$ where
$\beta \subseteq C_i$ with $r(C_i) \geq 2$. Also (as in Remark
\ref{top}), $f(E)=Eu$ where $u$ is a path in $C_i$ (the same
component of $G_{r-1}$ as $\beta$).
 By Proposition \ref{inv},
 there exists a path $\eta \subseteq C_i \subseteq G_{r-1}$ such
that $\bar{\eta} \beta \eta$ is a closed  Nielsen path at some
vertex $w'$ and that $\pi_1(f',w')$ has fixed subgroup of rank at
least $2$.
 If we now slide $E$ along $\eta $, we
get a new representative with an edge $E'$, such that $f'(E') = E'
[\bar{\eta} u f(\eta)]$. The new INP of height $r$ will be $E'
[\bar{\eta} \beta \eta] \bar{E'}$.

However, $$f(\beta) \simeq \bar{u} \beta u$$ and $$f( \bar{\eta}
\beta \eta ) \simeq \bar{\eta} \beta \eta .$$
Hence $\bar{\eta} u f(\eta) $ commutes with $\bar{\eta} \beta
\eta$ and since these are both closed paths we deduce that the
former is a power of the latter. (Note $\beta $ cannot be a proper
power.) Hence, $f'(E')= E' {\beta'}^k$ where $\beta'$ is a closed
Nielsen path at the vertex $w'$. Moreover $w'$ is fixed by $f'$,
$\pi_1(f',w')$ has fixed subgroup of rank at least $2$. Thus $f'$
is a map with the required properties and since no new vertices
were introduced, we may assume that the new representative has
minimal complexity. This concludes the proof of the claim.
\endproof

The idea now is to use the induction hypothesis on $C_1$ and
$C_2$. (Sliding an edge in $C_i$ is equivalent to sliding the same
edge in $G$). However it may be that $C_1,C_2$ contain valence one
vertices (these are not valence one in $G$ but in the $C_i$) so we
need to consider how this can arise. We proceed with a
representative of $\Phi$ which satisfies the conditions in the
claim above.

Consider first the case where $E$ is separating. Note that since
the INP of height $r$ is a closed path at $v$, this implies that
$\Sigma^w=\Sigma^w_{r-1}$ and hence that the fixed subgroup of
$\pi_1(f,w)$ is contained in $\pi_1(C_2,w)$. Since distinct INP's
start with distinct edges, we know that $w$ has valence at least
$2$ in $C_2$. Hence $C_2$ has no valence $1$ vertices and applying
the induction hypothesis we can assume that the Theorem holds for
every edge and vertex of $C_2$. (Note here that $f|_{C_2}$ is a
relative train track map and if it were not of minimal complexity
we could replace $f|_{C_2},C_2$ with some $f', C'$ via a homotopy
relative $w$. This is clearly not possible since it would
contradict the minimal complexity of $f,G$).

 To continue, if $r(C_1)=1$ we are done.
Also, if every vertex of $C_1$ has valence at least $2$, then we
are done since again we could apply the induction hypothesis to
$f|_{C_1},C_1$.

So there is only something to prove if $r(C_1) \geq 2$ and $C_1$
has a valence one vertex. Clearly, $v$ is the only vertex of $C_1$
which can have valence $1$. Let $e$ be the edge of $C_1$ whose
initial vertex is $v$. Since there is a closed INP at the vertex
$v$ contained in $C_1$, we deduce that the INP is of the form $e
\alpha \bar{e}$. Just as in the proof of the claim above, after
sliding $e$, we can assume that $f(e)=e \alpha^m$ for some $m$ and
that $\alpha $ is a closed Nielsen path.
 Note that the
terminal vertex of $e$ must have valence at least $3$, since
otherwise we could slide $e$ along an edge to produce a valence
$1$ vertex in $G$, contradicting the property of minimal
complexity. (This would also follow from the proof of the claim.)
Also, both endpoints of $e$ are fixed and that the corresponding
fixed subgroups have rank at least $2$ in $\pi_1(G)$.

Hence $f(C - \{e\}) \subseteq C -\{e\}$  and every vertex has
valence at least $2$. Thus we may apply our induction hypothesis
to $C - \{e \}$. The Theorem is then is then proved in this case,
since every edge of $G$ is either $E,e$ or in $C_1 \cup C_2$, and
every vertex is incident to one of these.

The same argument will apply when $E$ is non-separating, since we
may assume that $r(C_1) \geq 2$ as the case $n=2$ has already been
dealt with. \endproof

Let us call a representative of $\Phi$ which satisfies the
conclusions to Theorem \ref{goodrep} a good representative. One
immediate observation is that since every NP is closed, by Remark
\ref{vert}, the number of vertices of a good representative is
precisely $s(\Phi)$, the number of similarity classes with fixed
subgroup at least $2$. Thus if we start with an automorphism $\phi
\in Aut(F_n)$ which has fixed subgroup of rank $n$, then a good
representative of the outer automorphism induced by $\phi$ will
have exactly one vertex. Moreover, $\phi $ will be point
represented at that vertex and we recover the main Theorem of
\cite{CT2}. Methods used to analyze such automorphisms naturally
generalise to our situation. Hence the argument used in \cite{AM}
to solve the conjugacy problem for automorphisms with maximal
fixed subgroup can be applied with almost no changes to get:
\begin{Thm}
Given two outer automorphisms of maximal rank, $\Phi,\Psi \in
Out(F_n)$ in terms of images of a basis it is possible to decide
whether they are conjugate.
\end{Thm}

We shall show below, that any outer automorphism of maximal rank
can in fact be represented by a Dehn twist and the conjugacy
problem has been solved for these in \cite{CL}. Thus the only
advance made is an explicit algorithm when the automorphisms are
given in terms of images on a basis.

\section{Graphs of groups and Dehn twists}

We shall now show that any outer automorphism of $F_n$ of maximal
rank is represented by a Dehn Twist. We give a brief recap of the
objects involved, taken from \cite{CL}.

\begin{Dfn}
A graph of groups is given by $$\mathcal{G}=\{
\Gamma(\mathcal{G}),\{G_v \}_{e \in E(\mathcal{G})},\{f_e\}_{e \in
E(\mathcal{G})} \}$$
$\Gamma(\mathcal{G})$ is a finite connected graph,

$V(\mathcal{G})$ is the vertex set of $\Gamma(\mathcal{G})$,

$E(\mathcal{G})$ is the (oriented) edge set of $\Gamma(\mathcal{G})$,

$G_v$ is the vertex group at $v \in V(\mathcal{G})$,

$G_e$ is the edge group at $e \in E(\mathcal{G})$ and,

$m_e\co G_e \to G_{\tau(e)}$ is a monomorphism. (Here $\tau(e)$
denotes the terminal vertex of $e$ and hence $\tau{\bar{e}}$ the
initial vertex.)
\end{Dfn}

The {\em path group} $\Pi(\mathcal{G})$ is the free product of the
free group on the set $\{t_e : e \in E(\mathcal{G})\}$ with the
groups $G_v$, subject to the relations, 

(i)\qua $t_{\bar{e}}={t_e}^{-1}$ and

(ii)\qua $t_e m_e(a) {t_e}^{-1}=m_{\bar{e}}(a) \in G_{\tau{\bar{e}}}$
for all $a \in G_e, e \in E(\mathcal{G})$.

Every element of $\Pi(\mathcal{G})$ is given by a word $$W=r_0
t_1\ldots t_q r_q,$$ where each $t_i=t_{e_i}$ and $r_i$ is an
element of the free product of the $G_v$. Such a word is called a
loop at the vertex $v$ if $r_0, r_q \in G_v$,
$\tau(\bar{e_1})=\tau{e_q}=v$ and $\tau{e_i}=\tau{\bar{e_{i+1}}}$
with $r_i \in G_{\tau(e_i)}$ for all $i$ (taking subscripts modulo
$q$).

The set of loops at $v$ forms a subgroup of $\Pi(\mathcal{G})$
denoted, $\pi_1(\mathcal{G},v)$ and called the fundamental group
of $\mathcal{G}$ at $v$.

A Dehn twist, $\mathcal{D}$ on $\mathcal{G}$, with twistors $z_e$
is an automorphism of $\Pi(\mathcal{G})$ such that, $$
\mathcal{D}(t_e)= t_e m_e(z_e)$$ where $z_e$ is in the centre of
$G_e$ and $z_{\bar{e}}={z_e}^{-1}$. Extend $D$ to the whole of
$\Pi(\mathcal{G})$ by setting it equal to the identity on each
vertex group. Note that specifying the twistors $z_e$ is
sufficient to define $D$.

Since $D$ preserves incidence relations, it restricts to an
automorphism $D_v$ on $\pi_1(\mathcal{G})$ for any $v \in
V(\mathcal{G})$.

We shall say that an outer automorphism $\Phi \in Out(F_n)$ is
represented by a Dehn twist, $D$ on $\mathcal{G}$, at the vertex
$v$, if there is an isomorphism $\sigma: F_n \to
\pi_1(\mathcal{G},v)$ and a $\phi \in \Phi$ such that the
following diagram commutes, $$ \xymatrix{ F_n \ar[r]^{\sigma}
\ar[d]^{\phi} & \pi_1(\mathcal{G},v) \ar[d]^{D_v}
\\ F_n \ar[r]^{\sigma} & \pi_1(\mathcal{G},v).}$$
Note that $\pi_1(\mathcal{G},v)$ and $\pi_1(\mathcal{G},w)$ are
conjugate in $\Pi(\mathcal{G})$ and that under this ismorphism,
$D_v$ and $D_w$ define the same outer automorphism. Hence we may
refer to the outer automorphism induced by $D$. In particular, if
$\Phi $ is represented at the vertex $v$, then it will also be
represented at every other vertex, with different choices of $\phi
\in \Phi$.

 Now it is clear that if $\Phi $ is represented by a Dehn
twist $D$, then $\Phi $ and $D$ will have the same index, in fact
they will also have the same rank (as outer automorphisms). In
\cite{CL}, Corollary $7.7$ it is shown that every Dehn twist has
maximal index and no attracting fixed infinite words. Since this
is our definition of maximal rank, it follows that if $\Phi $ is
represented by a Dehn Twist then it has maximal rank. We now prove
the converse.

\begin{Thm}
Let $\Phi \in Out(F_n)$, $n \geq 2$ have maximal rank. Then $\Phi
$ is represented by a Dehn twist.
\end{Thm}

\proof Start with a relative train track map representative, $f,G$,
of $\Phi$ of minimal complexity satisfying the conclusion to
Theorem \ref{goodrep}. Form a graph $X$ from $G$ by deleting all
fixed edge loops. Thus $X$ is a connected subgraph of $G$ on the
same vertex set as $G$.

 We now define a graph of groups $\mathcal{G}$ with
graph $\Gamma(\mathcal{G})=X$. For each vertex $v$, let
$G_v=\fp{\pi_1(f,v)}$ and let each edge group $G_e$ be infinite
cyclic with generator $a_e$. Now each edge, $e$,  of
$\Gamma(\mathcal{G})=X$ is also an edge of $G$ so we can determine
the path $f(e)$ in $G$. By definition of $X$, $e$ is not a fixed
edge loop and so, by Proposition \ref{goodinp} and Theorem
\ref{goodrep}, up to orientation there is an INP $e \beta_e
\bar{e}$ and $f(e)=e \beta_e^{r_e}$ for some integer $r_e$ and a
$\beta_e$ which is a closed INP at the endpoint $\tau(e)$ of $e$.
We can then define a monomorphism $m_e\co  G_e \to G_{\tau(e)}$ by
mapping $a_e$ to $\beta_e$. Similarly, $f_{\bar{e}} \co
G_{\bar{e}} = G_e \to G_{\tau{\bar{e}}}$ will map $a_e$ to $e
\beta_e \bar{e}$. This completes the definition for $\mathcal{G}$.

Given an edge $e$ as above, set the twistor $z_{e}={a_e}^{r_e}$.
(Recall that $f(e)= e {\beta_e}^{r_e}$.) Since $G_e$ is abelian
this clearly lies in the centre. Then with
$\bar{z_e}={a_e}^{-r_e}$ we let $\mathcal{D}$ be the Dehn twist
based on these twistors.

It should be noted that if $e$ has height $k$ in $G$, with respect
to the stratification then there is a unique INP of height $k$ and
since $f,G$ has minimal complexity this INP will be closed and of
the form given by Proposition \ref{goodinp}. Thus the definition
of the maps $m_e$ for the graph of groups and the Dehn twist,
$\mathcal{D}$ is unambiguous.

Now let $\Pi(G)$ denote the free group on the (unoriented) edge
set of $G$ (so $e$ and $\bar{e}$ are considered inverse). Let
$\sigma $ map an edge $e$ of $G$ to the corresponding stable
letter $t_e$ of $\Pi(\mathcal{G})$ if $e$ is not a fixed edge
loop. If $e$ is a fixed edge loop at $v$, let $\sigma(e)$ be the
corresponding element of $G_v = \fp{\pi_1(f,v)}$. Since $\Pi(G)$
is free, $\sigma $ extends to a unique homomorphism from $\Pi(G)$
to $\Pi(\mathcal{G})$. We claim that this is an isomorphism.

One can easily define the inverse, $\sigma'$ from
$\Pi(\mathcal{G})$ to $\Pi(G)$. Simply let $\sigma'(t_e)=e$ for
every stable letter $t_e$ and for each $r \in G_v=
\fp{\pi_1(f,v)}$, let $\sigma'(r)$ be the unique product of INP's
representing it. We can extend $\sigma'$ to a well defined
homomorphism after an easy check to show that the relations in
$\Pi(\mathcal{G})$ are in the kernel of $\sigma'$. It is then a
matter of computation to show that the composition (in either
order)  of $\sigma $ and $\sigma'$ is the identity, by just
checking on generating sets. Thus $\sigma $ is an isomorphism as
claimed.

Now as $f$ maps edges to edge paths, it induces an endomorphism
$f_*$ on $\Pi(G)$. One can easily verify that the following
diagram commutes, $$\xymatrix{ \Pi(G) \ar[r]^{\sigma} \ar[d]^{f_*}
& \Pi(\mathcal{G}) \ar[d]^D \\ \Pi(G) \ar[r]^{\sigma} &
\Pi(\mathcal{G}).\\ }$$ The map $D$ is the Dehn twist defined
above. Note that $f_*$ is actually an automorphism of $\Pi(G)$ and
one can prove this either by induction on the number of strata in
$G$ or by observing the above diagram.

With this point of view, $\pi_1(G,v)$ is actually a subgroup of
$\Pi(G)$ and $\pi_1(f,v)$ is the restriction of $f*$ to this
subgroup. Since $\sigma$  preserves incidence relations one can
immediately deduce that $\sigma $ restricts to an isomorphism
between $\pi_1(G,v)$ and $\pi_1(\mathcal{G},v)$ and thus the above
commuting diagram restricts to, $$\xymatrix{ \pi_1(G,v)
\ar[r]^{\sigma} \ar[d]^{\pi_1(f,v)} & \pi_1(\mathcal{G},v)
\ar[d]^{D_v}
\\ \pi_1(G,v) \ar[r]^{\sigma} & \pi_1(\mathcal{G},v). }$$
This concludes the proof of the Theorem. \endproof

Note that in our proof we have done more than show a maximal rank
outer automorphism a represented by a Dehn twist. We have
represented the outer automorphism {\em when considered as a
groupoid automorphism} by a naturally equivalent Dehn twist, {\em
also considered as a groupoid automorphism}. This highlights the
strong connection between the two structures.

\section*{Acknowledgements}
This work was completed with the support of an EPSRC fellowship.

\Addresses\recd

\begin{thebibliography}


\bibitem{BFH} M. Bestvina and M. Feighn and M. Handel, The Tits
Alternative for $Out(F_n)$ I: Dynamics of Exponentially Growing
Automorphisms.  {\em Annals of Math.} 151 (2000), 517--623.

\bibitem{BH} M. Bestvina and M. Handel, Train tracks and automorphisms
of free groups. {\em Ann. of Math.} (2) 135 (1992), 1--51.

\bibitem{CL} M. M. Cohen and M. Lustig, The conjugacy problem for
Dehn twist automorphisms of free groups. {\em Comment Math. Helv.}
74 (1999) no. 2 179--200

\bibitem{CT1} D. J. Collins and E. Turner, An automorphism of a
free group of finite rank with maximal rank fixed point subgroup
fixes a primitive element. {\em J. Pure and Applied Algebra} 88
(1993) 43--49

\bibitem{CT2} D. J. Collins and E. Turner, All automorphisms of
free groups with maximal rank fixed subgroups. {\em Math. Proc.
Cambridge Phil. Soc.} 119 (1996), no. 4, 615--630


\bibitem{DS} J. L. Dyer and G. P. Scott, Periodic automorphisms of
free groups. {\em Comm. Alg.} 3 (1975), 195--201

\bibitem{GJLL} D. Gaboriau and A. Jaeger and G. Levitt and M. Lustig,
 An index for counting fixed points for automorphisms of free
 groups. {\em Duke Math. J.} 93 (1998), no. 3, 425--452.


\bibitem{GLL} D. Gaboriau and G. Levitt and M. Lustig, A
dendrological proof of the Scott conjecture for automorphisms of
free groups.  {\em Proc. Edinburgh Math. Soc.} 41 (1998), no. 2,
325--322

\bibitem{AM} A. Martino, Normal forms for automorphisms of maximal rank.
{\em Quart. J. Math.} 51 (2000), no. 4, 509--522



\end{thebibliography}
\end{document}